\documentclass[12pt]{article}
\usepackage{amssymb,amsmath,epsfig,graphics}

\newtheorem{corollary}{Corollary}
\newtheorem{example}{Example}

\newcommand{\qed}{{\hfill$\blacksquare$}}
\newcommand{\RR}{\mathbb R}
\newcommand{\CC}{\mathbb C}
\newcommand{\QQ}{\mathbb Q}
\newcommand{\ZZ}{\mathbb Z}

\newcommand{\NN}{\mathbb N}

\newcommand{\Ma}{\mathsf{M}}
\newcommand{\ma}{\mathsf{m}}

\newcommand{\CM}{\mathcal{M}}
\newcommand{\CQ}{\mathsf{Q}}
\newcommand{\cq}{\mathsf{q}}
\newcommand{\CK}{\mathsf{K}}

\newcommand{\tP}{\widetilde{P}}

\newcommand{\nr}{{$\sharp$ }}

\newcommand{\go}{{\omega }}

\newcommand{\gL}{{\Lambda }}

\newcommand{\af}{\hspace{-.8em}}
\newcommand{\dimer}[9]{\left[\hspace{-.5em}\begin{array}{ccccccccc}
i&\af h&\af g&\af f&\af e&\af d&\af c&\af b&\af a\\
#1&\af #2&\af #3&\af #4&\af #5&\af #6&\af #7&\af #8&\af #9
\end{array}\hspace{-.5em}\right]}

\begin{document}

\title{Mahler Measure, Eisenstein Series and Dimers.}

\author{Jan Stienstra\footnote{e-mail: stien{@}math.uu.nl}\\
\small{Mathematisch Instituut, Universiteit Utrecht, the Netherlands}\normalsize}
\date{}

\maketitle

\begin{abstract}
This note reveals a mysterious link between the partition function of certain dimer models on $2$-dimensional tori and the $L$-function of their spectral curves. It also relates the partition function in certain 
families of dimer models to Eisenstein series.
\end{abstract}

\section*{Introduction}

The \emph{logarithmic Mahler measure} $\ma (F)$ and the \emph{Mahler measure} $\Ma (F)$ of a Laurent polynomial $F(x,y)$ with complex coefficients are:
\begin{equation}\label{eq:mahler}
\ma (F):=\frac{1}{(2\pi i)^2}\oint\!\!\oint_{|x|=|y|=1}
\log |F(x,y)|\,
\frac{dx}{x}\frac{dy}{y}\,,\qquad
\Ma (F):=\exp (\ma (F))\;.
\end{equation}
Boyd \cite{Bo} gives a survey of many Laurent polynomials for which 
$\ma (F)$ equals 
(numerically to many decimal places) a `simple' non-zero rational number times the derivative at $0$ of the L-function of the projective plane curve  $Z_F$ defined by the vanishing of $F$:
\begin{equation}\label{eq:special L}
\ma(F)\cdot\QQ^*\;=\;\mathrm{L}'(Z_F,0)\cdot\QQ^*\,.
\end{equation}
Tables 2, 5, 1 in \cite{Bo} give the following families of cubic polynomials $\tilde{F}(X,Y,Z)$
for which (\ref{eq:special L}) was checked numerically with
$F(x,y)=(xy)^{-1}\tilde{F}(x,y,1)$ and with integer values of the parameter $s$,
\begin{eqnarray}
\label{eq:P6}    X^2Y+XY^2+X^2Z+XZ^2+Y^2Z+YZ^2-sXYZ&&\\
\label{eq:P3}     X^2Y+Y^2Z+Z^2X-sXYZ&&\\
\label{eq:P4}     X^2Y+XY^2+XZ^2+YZ^2-sXYZ&&
\end{eqnarray}
Deninger \cite{D} and Rodriguez Villegas \cite{RV} showed that the experimentally observed relations (\ref{eq:special L}) agree with predictions from the Bloch-Beilinson conjectures.
Rodriguez Villegas \cite{RV} gave a proof of (\ref{eq:special L}) for some values of $s$ in Examples (\ref{eq:P6}), (\ref{eq:P3}), (\ref{eq:P4}).

On the other hand, Kenyon, Okounkov and Sheffield
gave a formula for the \emph{partition function per fundamental domain} of a dimer model which is exactly the same as Formula 
(\ref{eq:mahler}) for the Mahler measure of the characteristic polynomial of that dimer model (see \cite{KOS} Theorem 3.5). 
In this paper we show examples of dimer models with characteristic polynomials 
(\ref{eq:P6}), (\ref{eq:P3}), (\ref{eq:P4}).

\emph{This suggests that there may be some mysterious link between the partition function of a dimer model and the $L$-function of its spectral curve.}

Our dimer models come in families with
the parameter $s$ explicitly related to the weights in the dimer model.
The results for the Mahler measures of the polynomials (\ref{eq:P6}),
(\ref{eq:P3}), (\ref{eq:P4}) presented in \cite{S1} now imply that, for $|s|$ sufficiently large, the partition functions in these families of dimer models  are
$\CQ_6^{-1}, \CQ_3^{-\frac{1}{3}}, \CQ_4^{-\frac{1}{2}}$ respectively, with
\begin{equation}\label{eq:prod463}
\begin{array}{lllrrr}
\CQ_6&=&
\displaystyle{\cq\prod_{n\geq 1}(1-\cq^n)^{(-1)^{n-1} n\chi_{-3}(n)}}\;,&
\CQ_3&=&
\displaystyle{\cq\prod_{n\geq 1}(1-\cq^n)^{9n\chi_{-3}(n)}}\;,\\[2ex]
\CQ_4&=&
\displaystyle{\cq\prod_{n\geq 1}(1-\cq^n)^{4n\chi_{-4}(n)}}\;,&&&
\end{array}
\end{equation}
where $\chi_{-3}(n)=0,\,1,\,-1$ if $ n\equiv 0,\,1,\,2\bmod 3$
and
$\chi_{-4}(n)=0,\,1,\,0,\,-1$ if $n\equiv 0,\,1,\,2,\,3\bmod 4$
and where $\cq$ is explicitly related to $s$.

The Eisenstein series in the title are the logarithmic derivatives of the products in (\ref{eq:prod463}):
$$
1+\sum_{n\geq 1}\chi_{-3}(n)\frac{n^2(-\cq)^n}{1-\cq^n}\,,\quad
1-9\sum_{n\geq 1}\chi_{-3}(n)\frac{n^2\cq^n}{1-\cq^n}\,,\quad
1-4\sum_{n\geq 1}\chi_{-4}(n)\frac{n^2\cq^n}{1-\cq^n}\,.
$$
The products and Eisenstein series for $\CQ_3$ and $\CQ_4$
were also studied by Ramanujan (see \cite{ABYZ}).
We find the similarity with McMahon's function
$$
M(\cq):=\prod_{n\geq 1}(1-\cq^n)^{-n}
$$
very intriguing, in particular because McMahon's function appears in
the partition functions in \cite{O, ORV}. 

\section{Dimers}
\begin{figure}
\begin{center}
\setlength\epsfxsize{7cm}
\epsfbox{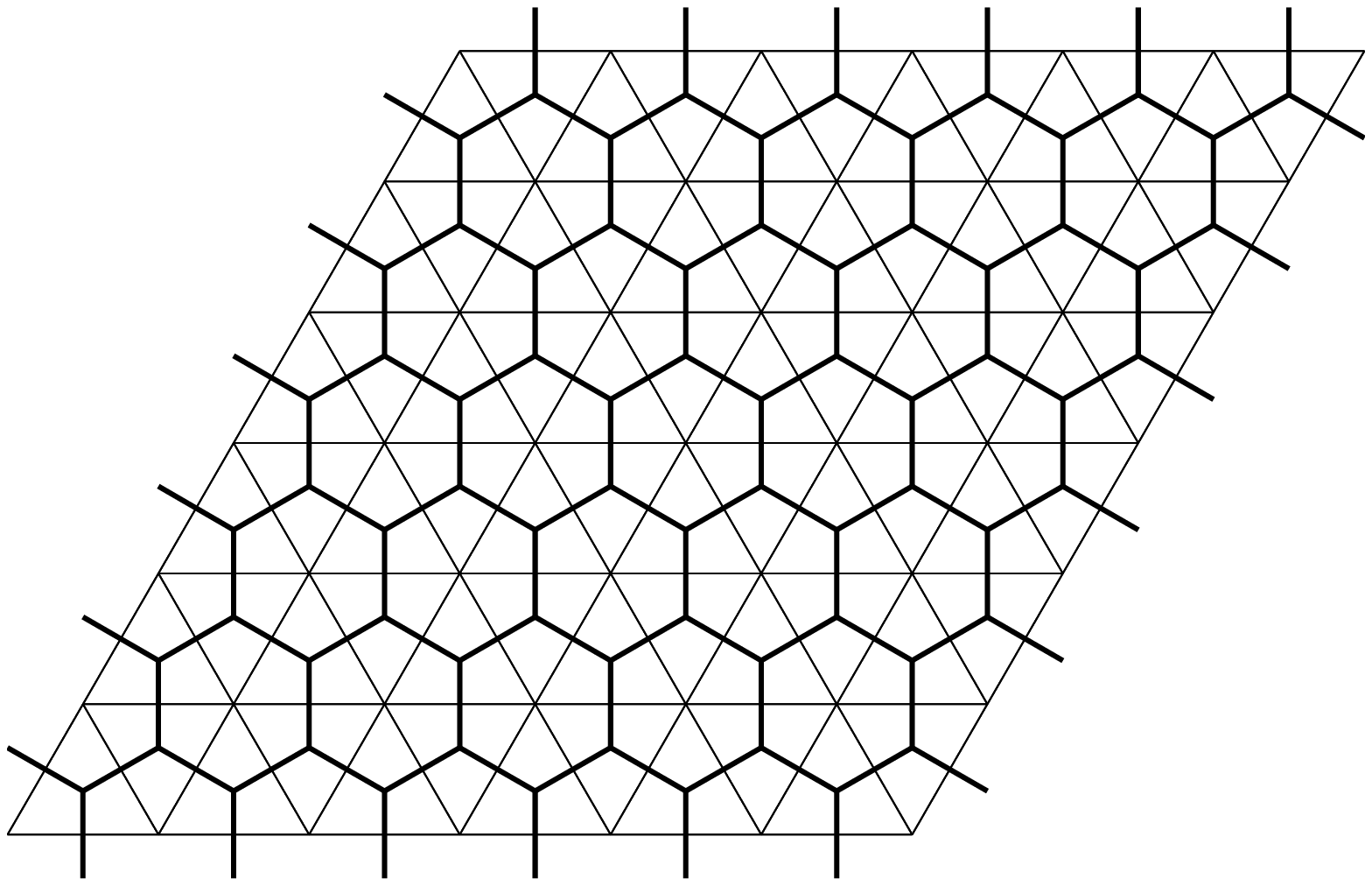}
\end{center}
\caption{\label{fig:honeycomb}}
\end{figure}

For general theory and background information on dimer models we refer to
\cite{KO, KOS, O, ORV}. We restrict ourselves to dimer models on the standard hexagonal graph $\Gamma$, which is the dual of the tessellation of the plane by regular triangles; see Figure \ref{fig:honeycomb} 

\begin{figure}
\hspace{-3em}
\setlength\epsfxsize{7cm}
\setlength\epsfysize{6cm}
\epsfbox{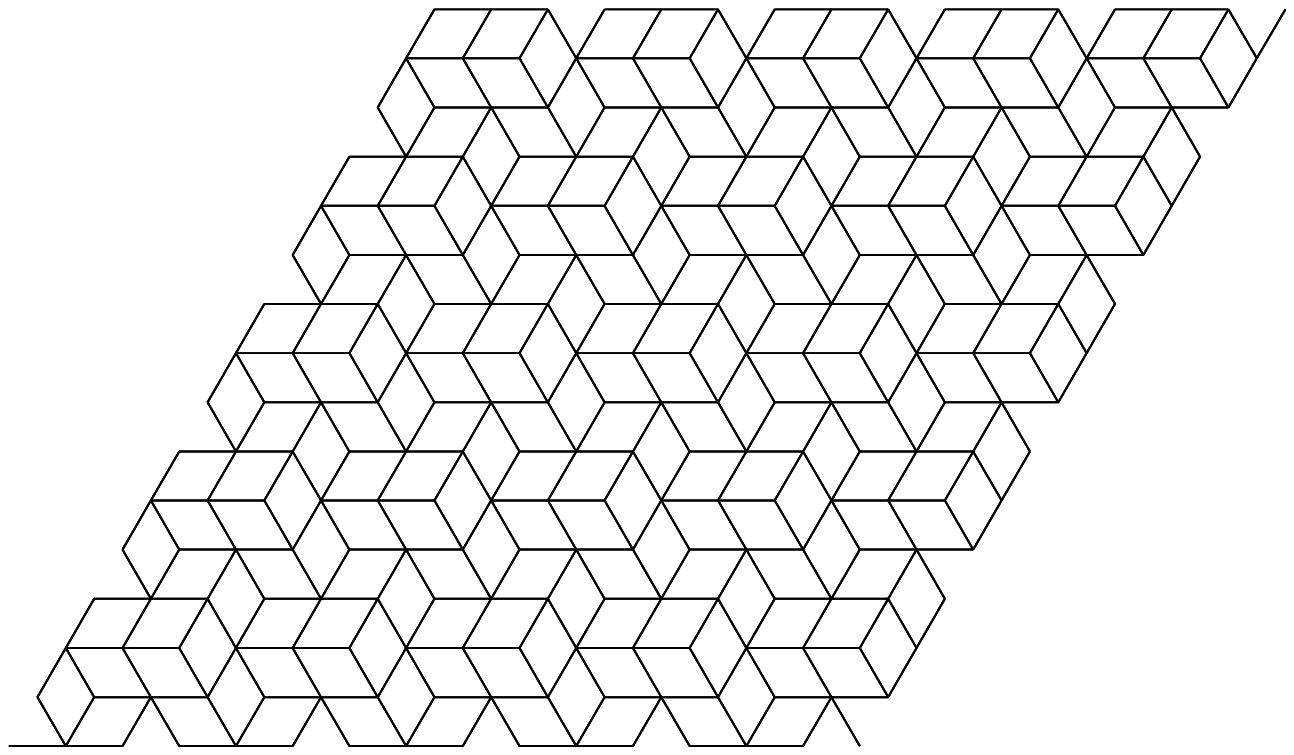}
\hspace{-3em}
\setlength\epsfxsize{7cm}
\setlength\epsfysize{6cm}
\epsfbox{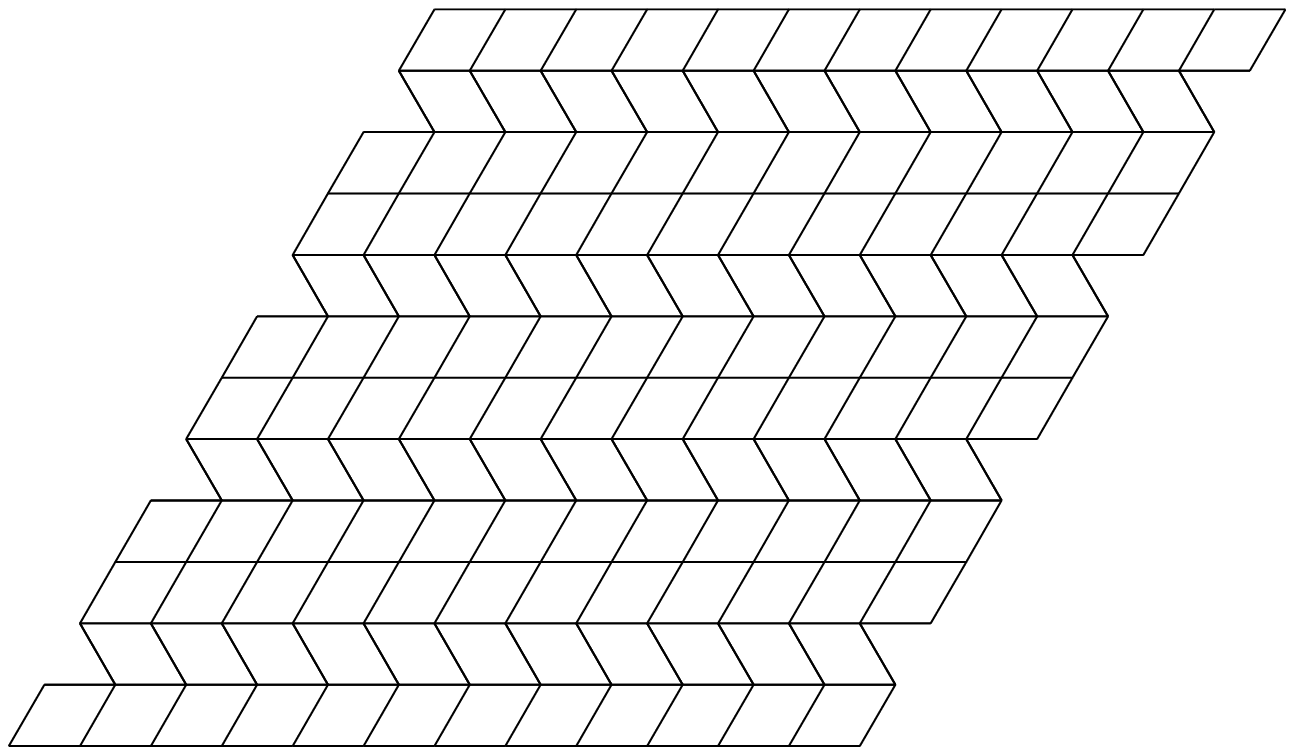}
\setlength\epsfxsize{7cm}
\setlength\epsfysize{6cm}
\epsfbox{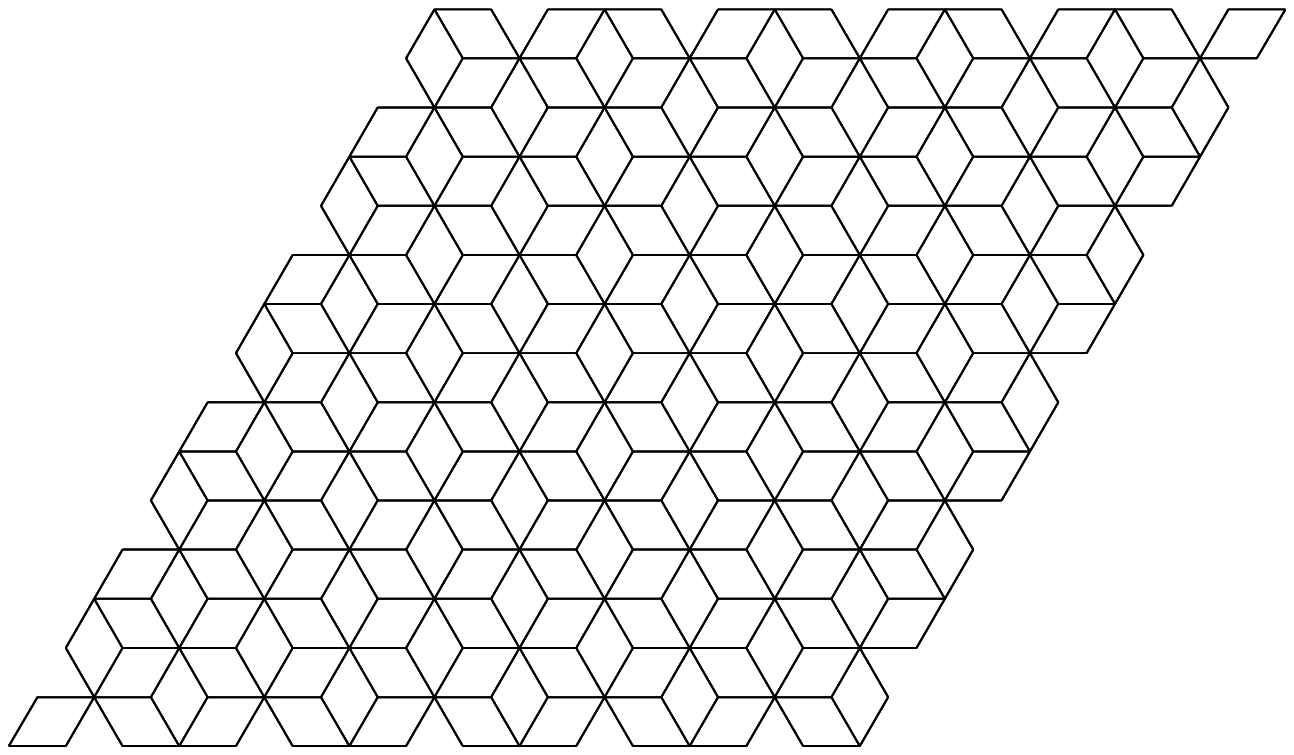}
\caption{\label{fig:tilings}
$180^\circ$ degree rotation of the top left picture
and $\pm 60^\circ$ degree rotations of the top right picture
give also legitimate tessellations.}
\end{figure}

A \emph{dimer configuration}, or perfect matching, on $\Gamma$
is a collection $\CM$ of its edges such that each vertex of $\Gamma$
is incident to exactly one edge from $\CM$. Each edge of $\Gamma$
intersects exactly one edge in the triangle tessellation and vice versa. Given a dimer configuration $\CM$ one can remove from the triangle tessellation all edges which intersect a $\Gamma$-edge belonging to $\CM$. The two triangles adjacent to such an edge glue together to a rhombus and as a result one gets a rhombus tessellation of the plane
(see Figure \ref{fig:tilings}).
Conversely every rhombus tessellation of the plane, in which each rhomb is the union of two adjacent triangles in the triangle tessellation, comes from a dimer configuration on $\Gamma$.

\begin{figure}
\begin{center}
\setlength\epsfxsize{5cm}
\epsfbox{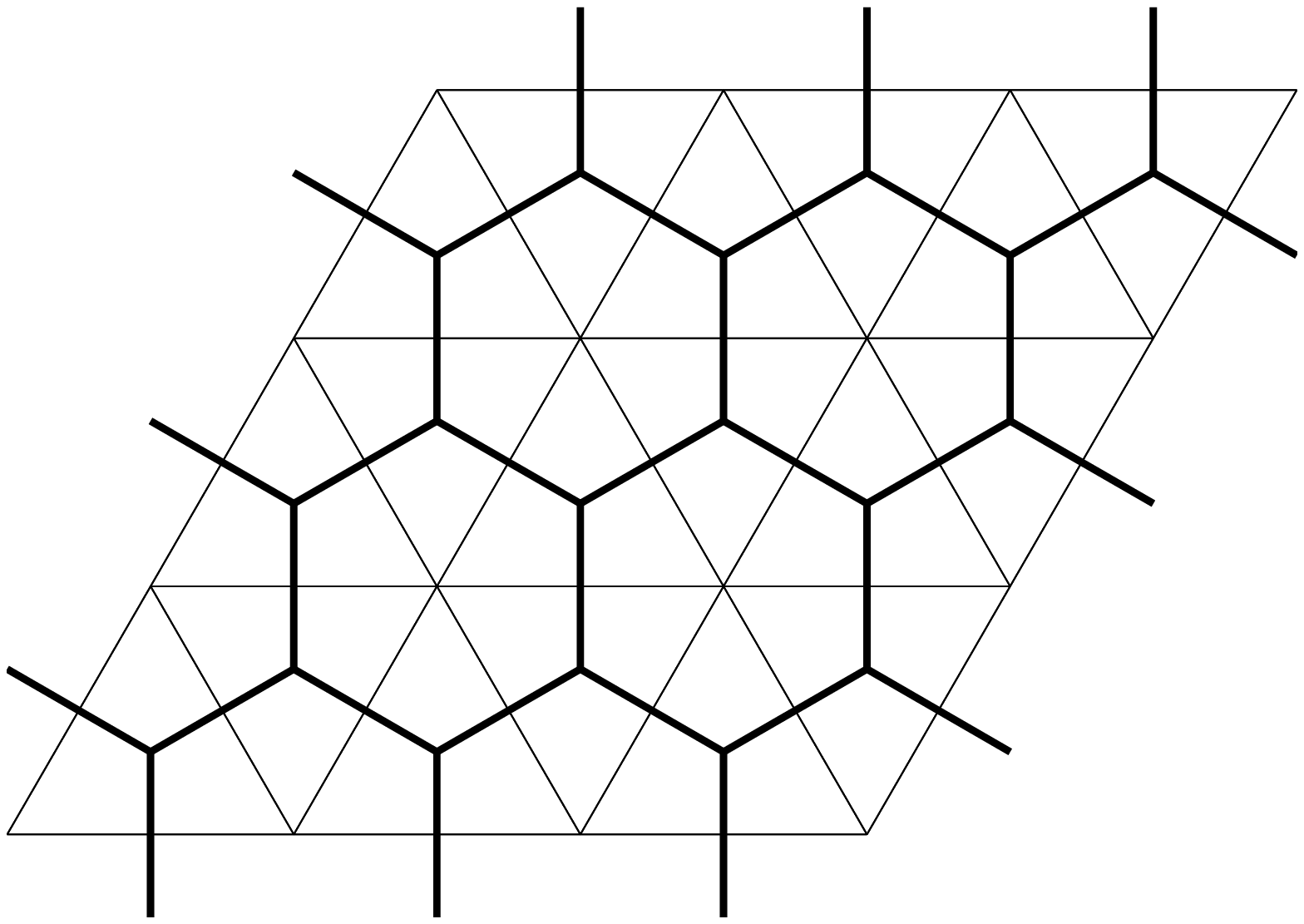}
\raisebox{9ex}{
$\qquad\begin{array}{cccccccc}
&&&a&&b&&c\\ &&A&&B&&C&\\[.8ex] &&d&&e&&f&\\ 
&D&&E&&F&&\\[.8ex] &g&&h&&i&&\\ G&&H&&I&&&
\end{array}$}

\footnotesize
$$
\CK:=
\left[ \begin{array}{ccccccccc}
Aa&Ba&0&0&0&0&Ga*x&0&0\\
0&Bb&Cb&0&0&0&0&Hb*x&0\\
Ac*y&0&Cc&0&0&0&0&0&Ic*x\\
Ad&0&0&Dd&Ed&0&0&0&0\\
0&Be&0&0&Ee&Fe&0&0&0\\
0&0&Cf&Df*y&0&Ff&0&0&0\\
0&0&0&Dg&0&0&Gg&Hg&0\\
0&0&0&0&Eh&0&0&Hh&Ih\\
0&0&0&0&0&Fi&Gi*y&0&Ii
\end{array}\right]\,.
$$
\end{center}
\caption{\label{fig:fundamental}
\footnotesize{
The matrix entries give an abstract, tautological weighting of the graph edges.
The $*x$ (resp. $*y$) marks the graph edges for which the top-bottom
(resp. left-right) sides of the fundamental parallellogram must be identified.}
\normalsize}
\end{figure}

View the plane as the complex plane $\CC$ and let
$\go:=e^{2\pi i/3}$. Then the vertices of the triangulation
are precisely the elements of the lattice
$$
\gL:=\{ a+b\go\,|\, a,b\in\ZZ\}\,.
$$ 
We first look at dimer configurations which are invariant under translations by vectors from the lattice $3\gL$. One may also say that these are dimer configurations on the graph $\Gamma/3\gL$, embedded in the torus $\CC/3\gL$. Figure \ref{fig:fundamental} shows 
the fundamental domain
and a labeling scheme for the vertices of the graph c.q. the triangles of the tessellation. It also shows
the same adjacency and labeling structure in the matrix $\CK$,
which is the matrix of the \emph{Kasteleijn operator} of the dimer model
on $\Gamma/3\gL$. The determinant 
\begin{equation}\label{eq:charpol}
P(x,y):=\det\CK
\end{equation}
is called the \emph{characteristic polynomial} and
the zero locus of $P(x,y)$ in $\CC^{*2}$ is called the \emph{spectral curve} of the dimer model (see \cite{KO} \S 2.1).
A straightforward computation gives the characteristic polynomial shown in Figure \ref{fig:charpol}. 

Next we want to specify a $3\gL$-invariant weight function $W$ on the edges of $\Gamma$ so that the characteristic polynomial
of the dimer model with these weights is a constant times one of the  polynomials (\ref{eq:P6}), (\ref{eq:P3}), (\ref{eq:P4}).

\begin{figure}
\footnotesize
$$
\begin{array}{l}
\tP(X,Y,Z):=Z^3P(\textstyle{\frac{X}{Z},\frac{Y}{Z}})=\\[3ex]
\hspace{-.5em}
X^3\dimer{F}{E}{D}{C}{B}{A}{I}{H}{G}
+Y^3\dimer{G}{I}{H}{D}{F}{E}{A}{C}{B}
+Z^3\dimer{I}{H}{G}{F}{E}{D}{C}{B}{A}\\[3ex]
+XZ^2( \dimer{F}{H}{G}{C}{E}{D}{I}{B}{A}
    + \dimer{I}{E}{G}{F}{B}{D}{C}{H}{A}  
    + \dimer{I}{H}{D}{F}{E}{A}{C}{B}{G})\\[3ex]
+X^2Z( \dimer{F}{E}{G}{C}{B}{D}{I}{H}{A} 
     + \dimer{F}{H}{D}{C}{E}{A}{I}{B}{G} 
     + \dimer{I}{E}{D}{F}{B}{A}{C}{H}{G})\\[3ex]
+YZ^2( \dimer{G}{I}{H}{F}{E}{D}{C}{B}{A}
     + \dimer{I}{H}{G}{D}{F}{E}{C}{B}{A} 
     + \dimer{I}{H}{G}{F}{E}{D}{A}{C}{B})\\[3ex]
+Y^2Z( \dimer{G}{I}{H}{D}{F}{E}{C}{B}{A}
     + \dimer{G}{I}{H}{F}{E}{D}{A}{C}{B}
     + \dimer{I}{H}{G}{D}{F}{E}{A}{C}{B})\\[3ex] 
+X^2Y( \dimer{G}{E}{D}{C}{F}{A}{I}{H}{B}
    + \dimer{F}{E}{H}{D}{B}{A}{I}{C}{G}
    + \dimer{F}{I}{D}{C}{B}{E}{A}{H}{G})\\[3ex]
+XY^2( \dimer{G}{E}{H}{D}{F}{A}{I}{C}{B}
    + \dimer{G}{I}{D}{C}{F}{E}{A}{H}{B}
    + \dimer{F}{I}{H}{D}{B}{E}{A}{C}{G})\\[3ex]
-XYZ(\;\dimer{G}{E}{H}{C}{F}{D}{I}{B}{A}
+ \dimer{G}{H}{D}{C}{F}{E}{I}{B}{A} 
+ \dimer{G}{E}{H}{F}{B}{D}{I}{C}{A}
\\[3ex] \hspace{3em}  +\dimer{F}{H}{G}{D}{B}{E}{I}{C}{A}
+ \dimer{G}{H}{D}{F}{B}{E}{I}{C}{A}
+ \dimer{F}{I}{G}{D}{B}{E}{C}{H}{A} 
\\[3ex] \hspace{3em} + \dimer{G}{I}{D}{F}{B}{E}{C}{H}{A}
+ \dimer{F}{H}{G}{D}{E}{A}{I}{C}{B} 
+ \dimer{G}{H}{D}{F}{E}{A}{I}{C}{B} 
\\[3ex] \hspace{3em} + \dimer{F}{I}{G}{C}{E}{D}{A}{H}{B} 
+ \dimer{I}{E}{G}{C}{F}{D}{A}{H}{B} 
+ \dimer{F}{I}{G}{D}{E}{A}{C}{H}{B} 
\\[3ex] \hspace{3em} + \dimer{G}{I}{D}{F}{E}{A}{C}{H}{B} 
+ \dimer{I}{E}{G}{D}{F}{A}{C}{H}{B} 
+ \dimer{F}{I}{H}{C}{E}{D}{A}{B}{G} 
\\[3ex] \hspace{3em} + \dimer{I}{E}{H}{C}{F}{D}{A}{B}{G} 
+ \dimer{I}{H}{D}{C}{F}{E}{A}{B}{G} 
+ \dimer{F}{I}{H}{D}{E}{A}{C}{B}{G}
\\[3ex]  \hspace{3em} + \dimer{I}{E}{H}{D}{F}{A}{C}{B}{G} 
+ \dimer{I}{E}{H}{F}{B}{D}{A}{C}{G} 
+ \dimer{I}{H}{D}{F}{B}{E}{A}{C}{G}).
\end{array}
$$
\caption{\label{fig:charpol}}
here
$\dimer{F}{E}{D}{C}{B}{A}{I}{H}{G}$ denotes the product
$Fi*Eh*Dg*Cf*Be*Ad*Ic*Hb*Ga$ of the weights.
It also shows the dimer configuration as a map from
$\nabla$- to $\Delta$-triangles.
\normalsize
\end{figure}
So, the coefficients of certain monomials must be $0$ and the coefficients of all other monomials different from $XYZ$ must be equal and non-zero.
Moreover, the ratio between the coefficient of $XYZ$ and the other non-zero coefficient should be $-s$. Altogether these conditions constitute $9$ equations on the weights.
Which $9$ equations exactly depends on the polynomial one wants to reconstruct.

The graph $\Gamma/3\gL$ has $27$ edges and $18$ vertices.
If one multiplies the weights of the three edges incident to a vertex by the same non-zero number $k$, the whole characteristic polynomial gets multiplied by $k$. So there is a torus ${\CC^*}^{18}$ acting on the solutions of the $9$ equations, but the $1$-dimensional subtorus
given by multiplication factors $k$ for the nine $\Delta$-vertices and 
$k^{-1}$ for the nine $\nabla$-vertices ($k\in\CC^*$) acts trivially. 
So, for fixed $s$, the solution space of the $9$ equations modulo the action of ${\CC^*}^{18}$ has dimension at least
$1=27-9-17$.

Without further delving into general solutions we now exhibit just three examples from a larger collection we found by trial and error. Correctness of these examples can easily be checked by direct inspection of the general characteristic polynomial $\tP(X,Y,Z)$ in
Figure \ref{fig:charpol}.

\begin{figure}
\begin{center}
\setlength\epsfxsize{10cm}
\epsfbox{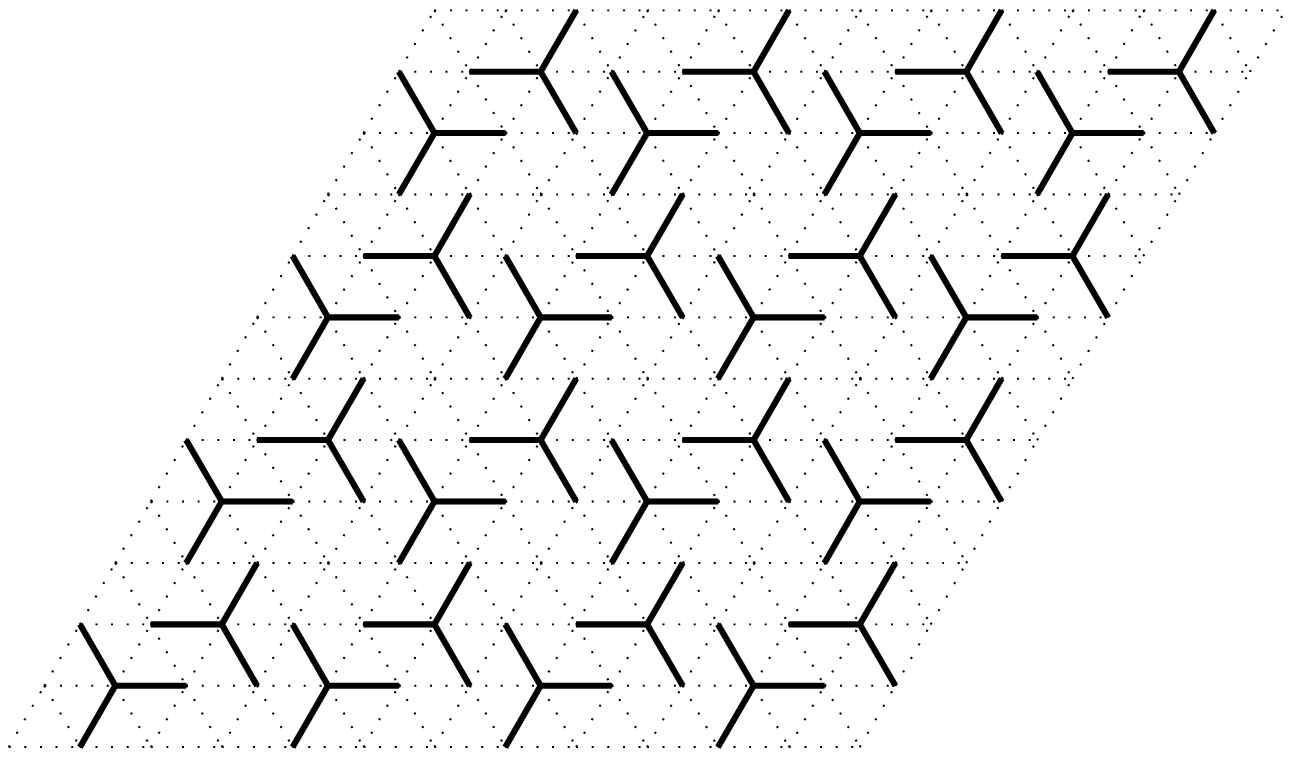}
\end{center}
\caption{\label{fig:obstacles}}
\end{figure}

As a common structure in these examples there are only four edge weights ($0,\,1,\,w,\,m$) and $s$ is a function of $m$. The parameter $w$ confirms that for fixed $s$ the solution space modulo ${\CC^*}^{18}$ has dimension $\geq 1$. Edges with weight $0$ can not occur in a dimer configuration. Thus if the dimer configuration is presented as a rhombus tessellation
the $0$-weight edges of $\Gamma$ turn into edges of the tessellation
which \emph{must be present}; it is like boundary conditions. This is shown for Example \ref{exa3} in
Figure \ref{fig:obstacles}. Any rhombus tessellation which respects the boundary conditions is allowed.
Figure \ref{fig:tilings} shows the $3\gL$-periodic dimer configurations
for Example \ref{exa3}.

\begin{example}\label{exa6}
Setting $W(Hg)=W(Dd)=W(Eh)=0$, $\;W(Hb)=W(Ee)=W(Df)=m$, $\;W(Ii)=W(Ic)=W(Gi)=W(Aa)=W(Ac)=W(Ga)=1$ and all remaining weights $=w$
yields a characteristic polynomial of type (\ref{eq:P6}):
$$
mw^6(X^2Y+XY^2+X^2Z+XZ^2+Y^2Z+YZ^2)-(4 + 3m + 3m^2 + m^3)w^6XYZ
$$
with $s=(4 + 3m + 3m^2 + m^3)m^{-1}$.
\end{example}

\begin{example}\label{exa3}
Setting $W(Hg)=W(Dd)=W(Eh)=W(Be)=W(Cb)=W(Ff)=0$,
$\;W(Hb)=W(Ee)=W(Df)=m$,
$\;W(Ii)=W(Ic)=W(Gi)=W(Aa)=W(Ac)=W(Ga)=1$ and all remaining weights $=w$
yields a characteristic polynomial of type (\ref{eq:P3}):
$$
mw^6(X^2Z+XY^2+YZ^2)-(2+m^3)w^6XYZ
$$
with $s=(2+m^3)m^{-1}$.
\end{example}

\begin{example}\label{exa4}
Setting $W(Aa)=W(Bb)=W(Cf)=W(Gi)=W(Fe)=0$, $\;W(Cb)=W(Ff)=W(Be)=1$, $\;W(Cc)=W(Fi)=W(Ba)=m$ and all remaining weights $=w$
yields a characteristic polynomial of type (\ref{eq:P4}):
$$
mw^6(X^2Y+XY^2+XZ^2+YZ^2)-(2 + m^2+m^3)w^6XYZ
$$
with $s=(2 + m^2+m^3)m^{-1}$.
\end{example}
\textbf{Remark:} Among the examples in \cite{KOS} there is a dimer model on the square-octagon graph with characteristic polynomial of type (\ref{eq:P4}):
$x+x^{-1}+y+y^{-1}-s$.

\

As in \cite{KOS} \S3.2, we now fix a $3\gL$-invariant $\RR$-valued
weight function $W$ on the edges of $\Gamma$ and consider for a positive integer $n$ the dimer configurations on $(\Gamma,W)$ which are invariant under translations by vectors in the lattice $3n\gL$. So, the preceding discussion and examples
concern the case $n=1$. Theorem 3.3 of \cite{KOS} expresses the characteristic polynomial $P_n(x,y)$ of the dimer model on $(\Gamma/3n\gL,W)$ in terms of the characteristic polynomial $P(x,y)$ of the dimer model on $(\Gamma/3\gL,W)$:
$$
P_n(x,y)\,=\,\prod_{u^n=x}\prod_{v^n=y} P(u,v)\,.
$$
Corollary 3.4 of \cite{KOS} uses this to compute the \emph{partition function} $Z(\Gamma/3n\gL,W)$ for the dimer model on $(\Gamma/3n\gL,W)$:
\begin{eqnarray}
\nonumber
Z(\Gamma/3n\gL,W)&=&
\frac{1}{2}(-Z_n^{(0,0)}+Z_n^{(0,1)}+Z_n^{(1,0)}+Z_n^{(1,1)})\\[2ex]
\hspace{-5em}\textrm{where}\hspace{5em}
\label{eq:partition n}
Z_n^{(a,b)}&=&\prod_{u^n=(-1)^a}\prod_{v^n=(-1)^b} P(u,v)\,.
\end{eqnarray}
Finally, \cite{KOS} Thm. 3.5 computes the \emph{partition function per fundamental domain} $Z$ for periodic dimer models on $(\Gamma,W)$
with period lattice $3n\gL$, $n\in\NN$:
\begin{eqnarray}
\nonumber
\log Z&:=&\lim_{n\to\infty}\frac{1}{n^2}\log Z(\Gamma/3n\gL,W)\\
\label{eq:partition function} &=&
\frac{1}{(2\pi i)^2}\oint\!\!\oint_{|x|=|y|=1}
\log |P(x,y)|\,
\frac{dx}{x}\frac{dy}{y}\quad=:\quad\ma(P)\,.
\end{eqnarray}
\emph{Note}: $|P(x,y)|=P(x,y)$ in (\ref{eq:partition function}), because the polynomial $P(x,y)$ has real coefficients and (\ref{eq:partition n}) contains with every factor also its complex conjugate.

\

\begin{corollary}
The partition function per fundamental domain is equal to the Mahler measure of the characteristic polynomial:
$$
Z=\Ma(P)\,.
$$
\qed\end{corollary}

We now combine the above examples with Examples \nr 6, \nr 3, \nr 4 of \cite{S1}.

\setcounter{example}{0}

\begin{example}\textup{\textbf{(continued)}}
According to \cite{S1} Example \nr 6 the Mahler measure of the polynomial
$$
(X+Y)(Y+Z)(Z+X)-tXYZ
$$
with $|t|$ sufficiently large, is $\CQ_6^{-1}$ where
$$
\CQ_6=
\displaystyle{\cq\prod_{n\geq 1}(1-\cq^n)^{(-1)^{n-1} n\chi_{-3}(n)}}
\quad\textrm{and}\quad
t=\frac{\eta(\cq^2)^3\eta(\cq^3)^9}{\eta(\cq)^3\eta(\cq^6)^9}
$$
and $\chi_{-3}(n)=0,\,1,\,-1$ if $ n\equiv 0,\,1,\,2\bmod 3$.
Thus, if we rescale the dimer model in Example \ref{exa6} by multiplying all edge-weights by $(mw^6)^{-\frac{1}{9}}$, the partition function per fundamental domain for the rescaled dimer model is $\CQ_6^{-1}$ with
$$
t-2=s=(4 + 3m + 3m^2 + m^3)m^{-1}\,.
$$
\end{example}

\begin{example}\textup{\textbf{(continued)}}
According to \cite{S1} Example \nr 3 the Mahler measure of the polynomial
$$
X^2Y+Y^2Z+Z^2X-tXYZ
$$
with $|t|$ sufficiently large, is $\CQ_3^{-\frac{1}{3}}$ where
$$
\CQ_3=
\displaystyle{\cq\prod_{n\geq 1}(1-\cq^n)^{9n\chi_{-3}(n)}}
\quad\textrm{and}\quad
t^3=27+\frac{\eta(\cq)^{12}}{\eta(\cq^3)^{12}}\,.
$$
Thus, if we rescale the dimer model in Example \ref{exa3} by multiplying all edge-weights by $(mw^6)^{-\frac{1}{9}}$, the partition function per fundamental domain for the rescaled dimer model is $\CQ_3^{-\frac{1}{3}}$ with
$$
t=s=(2+m^3)m^{-1}\,.
$$
\end{example}

\begin{example}\textup{\textbf{(continued)}}
According to \cite{S1} Example \nr 4 the Mahler measure of the polynomial
$$
(X+Y)(XY+Z^2)-tXYZ
$$
with $|t|$ sufficiently large, is $\CQ_4^{-\frac{1}{2}}$ where
$$
\CQ_4=
\displaystyle{\cq\prod_{n\geq 1}(1-\cq^n)^{4n\chi_{-4}(n)}}
\quad\textrm{and}\quad
t=\frac{\eta(\cq^2)^{12}}{\eta(\cq)^4\eta(\cq^4)^8}
$$
and
$\chi_{-4}(n)=0,\,1,\,0,\,-1$ if $n\equiv 0,\,1,\,2,\,3\bmod 4$.
Thus, if we rescale the dimer model in Example \ref{exa4} by multiplying all edge-weights by $(mw^6)^{-\frac{1}{9}}$, the partition function per fundamental domain for the rescaled dimer model is $\CQ_4^{-\frac{1}{2}}$ with
$$
t=s=(2 + m^2+m^3)m^{-1}\,.
$$
\end{example}

%



\begin{thebibliography}{99}
\bibitem{ABYZ}
Ahlgren, S., B. Berndt, A. Yee, A. Zaharescu,
\textit{Integrals of Eisenstein series and derivatives of L-functions},
International Math. Research Notices 32 (2002) 1723--1738.
\bibitem{Bo}
Boyd, D., \textit{Mahler's measure and special values of L-functions},
Experimental Math. vol. 7 (1998) 37--82
\bibitem{D}
Deninger, D. \textit{Deligne periods of mixed motives, K-theory and the entropy of certain $\ZZ^n$-actions}, J. Amer. Math. Soc. 10 (1997) 259--281
\bibitem{KO}
Kenyon, R., A. Okounkov, \textit{Planar dimers and Harnack curves},
arXiv:math.AG/0311062
\bibitem{KOS}
Kenyon, R., A. Okounkov, S. Sheffield, \textit{Dimers and Amoebae},
arXiv:math-ph/0311005
\bibitem{O}
Okounkov, A., \textit{Random surfaces enumerating algebraic curves},
arXiv:math-ph/0412008
\bibitem{ORV}
Okounkov, A., N. Reshetikin, C. Vafa, \textit{Quantum Calabi-Yau and Classical Crystals},
arXiv:hep-th/0309208
\bibitem{RV}
Rodriguez Villegas, F., \textit{Modular Mahler measures I}, 
 Topics in number theory (University Park, PA, 1997), Ahlgren, S., G. Andrews, K. Ono (eds) 17--48,
Math. Appl., 467,
Kluwer Acad. Publ., Dordrecht, 1999. See 
also: http://www.ma.utexas.edu/users/villegas/research.html
\bibitem{S1}
Stienstra, J., \textit{Mahler Measure Variations, Eisenstein Series and Instanton Expansions}, these proceedings
\end{thebibliography}
\end{document}